\documentclass[preprint,11pt]{elsarticle}

\usepackage{amssymb}
\usepackage{graphicx}
\usepackage{amsfonts}
\usepackage{hyperref,amsmath,amsthm}%lineno,
\usepackage{color}

\theoremstyle{definition}

\newtheorem{theorem}{Theorem}[section]
\newtheorem{lemma}{Lemma}[section]

\theoremstyle{remark}

\journal{Elsevier}%Applied and Computational Harmonic Analysis

\def\bbe{{\Bbb{E}}}

\begin{document}

\begin{frontmatter}

%% Title, authors and addresses

%% use the tnoteref command within \title for footnotes;
%% use the tnotetext command for theassociated footnote;
%% use the fnref command within \author or \address for footnotes;
%% use the fntext command for theassociated footnote;
%% use the corref command within \author for corresponding author footnotes;
%% use the cortext command for theassociated footnote;
%% use the ead command for the email address,
%% and the form \ead[url] for the home page:
%% \title{Title\tnoteref{label1}}
%% \tnotetext[label1]{}
%% \author{Name\corref{cor1}\fnref{label2}}
%% \ead{email address}
%% \ead[url]{home page}
%% \fntext[label2]{}
%% \cortext[cor1]{}
%% \address{Address\fnref{label3}}
%% \fntext[label3]{}

%\title{Forecast Oil Market Share under the COVID-19 Pandemic Using Two-stage Stochastic Equilibrium}

\title{Equilibrium Oil Market Share under the COVID-19 Pandemic}
%% use optional labels to link authors explicitly to addresses:
%% \author[label1,label2]{}
%% \address[label1]{}
%% \address[label2]{}

%\author{}
%
%\address{}
\author[1]{Xiaojun Chen\corref{mycorrespondingauthor}}
%\ead[url]{www.elsevier.com}
\cortext[mycorrespondingauthor]{Corresponding author}
\ead{maxjchen@polyu.edu.hk}

\author[2]{Yun Shi}
\ead{admin@bbinnovative.com}

\author[1]{Xiaozhou Wang}%\textsuperscript{a}}
\ead{xzhou.wang@connect.polyu.hk}

\address[1]{Department
of Applied Mathematics, The Hong Kong Polytechnic University, Hong Kong, China}
\address[2]{ Blue Balloon Innovative Limited,
%16A-B,
Charmhill Center,
% 50 Hillwood Road,
Tsimshatsui, Hong Kong, China}

\begin{abstract}
Equilibrium models for energy markets under uncertain demand and supply have attracted considerable attentions. This paper focuses on modelling crude oil market share under the COVID-19 pandemic using two-stage stochastic equilibrium.
We describe  the uncertainties in the demand and supply
  by random variables  and provide two types of production decisions (here-and-now and wait-and-see). The here-and-now decision in the first stage does not depend on the outcome of random events to be revealed in the future and the wait-and-see decision in the second stage is allowed to depend on the random events in the future and adjust the feasibility of the here-and-now  decision in rare unexpected scenarios such as those observed during the  COVID-19 pandemic. We develop a fast algorithm to find a solution of the two-stage stochastic equilibrium. We show the robustness of the two-stage stochastic equilibrium model for forecasting the oil market share using the real market data from January 2019 to May 2020.
\end{abstract}


\end{frontmatter}

%\linenumbers

\section{Introduction}
In this paper, we consider the two-stage stochastic equilibrium model under uncertain demand and supply for forecasting the market share  in an oligopoly market. This model is developed based on noncooperative game theory which allows us to investigate  how the market share of one  product (oil, electricity, steel, etc.) depends on the strategies of a few agents whose decisions against and affect each other. A Nash-equilibrium for oligopolies is expected, in which costs and benefits are balanced so that no agent can gain by changing only their own strategies. The uncertainties in the future events are represented by random variables in the model. The COVID-19 pandemic has  significant impact  on the global energy industry  and arouse great challenging in maintaining market stability. Whether we can expect a Nash-equilibrium for crude oil market share during the COVID-19 pandemic is an interesting question.

In traditional  oil industry, each producer has normally two strategies for maximizing its profits. One is to limit supply which leads to high oil price and this approach allows high-cost producers to remain profitable. The other is to drive up production and squeeze out high-cost producers, e.g., US shale oil, under expected low price market environment. It has been shown that under certain conditions, one strategy can be more appropriate than the other \cite{Behar}.
However, the change in production  had inconsistent market responses due to the complex and uncertain nature of oil market \cite{Ji,Panda}.
Furthermore, it has been observed that even during the periods of violently volatile in oil market, the relative market share remained stable. See Tables 3-4 for monthly oil market share from January 2019 to May 2020.

The industry of oil is no stranger to volatilities in price, supply and demand \cite{Juvenal,Yao}.
Economic growth can be largely affected by oil volatility, and certain countries may be more sensitive than the others \cite{Kristie,Nordhaus,Van}.
Most of these changes, the so-called oil shocks, were caused by the occurrence of major events of global economics and dramatic trends shifting that affects supply or/and demand sides in the crude oil market, and the situation was no different in the latest oil shock during COVID-19 pandemic. Events of the crude oil market are capable of delivering significant impact on energy markets as well as non-energy markets
\cite{Ji2012,Kan,Narayan}.
It is believed, however, despite different causes of oil shock, their consequences for economy are very similar \cite{Hamilton}.
The most apparent observation on oil uncertainty is its price, and it is naturally assumed as of certain importance to market participants.
However, being an oligopoly market implies that the price is largely manipulated by major producers. In particular, these producers may choose to sell greater quantity of oil at a lower price over the strategies of limiting supply to boost price, because the former can be more profitable in long run. There are many attempts on predicting crude oil price recently using tools of machine learning and artificial intelligence \cite{BarunA,Chiroma,Wang}.

%There exist different standards on the quality of traded oil, hence multiple benchmark prices of oil co-exist. It is worth noting that some of the benchmark prices, e.g., the West Texas Intermediate (WTI), are to be settled in future markets hence representing anticipated values in the future.
%However, the trends of majorities of such benchmarks are highly correlated and make little analytical differences.

To analyze the oil market share under uncertain supply and demand,  we model the production decisions as a solution of a two-stage stochastic game.
The two-stage stochastic game  is suitable to reflect the complexity of the trading process while market demand needs to be predicted/estimated at the second stage to guide the production at the first stage for a certain period of time in the future.
In the two-stage stochastic game, we have two types of decisions: ``here-and-now''
and ``wait-and-see'' for each agent.   Endowed with market data of production, supply, demand, cost and price of the commodity, each agent chooses its  here-and-now decision and wait-and-see decision with and without observation of uncertainties in the future events, respectively. For each agent, the objective of the first stage  is to maximize
its expected utilities under capacity constraints, and the objective of the second stage  is to maximize  its utilities subject to constraints endowed with its here-and-now decision and recourse term for almost every realization. The structure of the two-stage stochastic game is a natural way to represent
the equilibrium state under uncertain environment.  Each agent must make a here-and-now decision for a production plan before knowing the future events. The here-and-now decision will affect future revenues, cost and feasibilities over the decision horizon. The wait-and-see decision is made by taking the uncertain demand, supply and price  in the future market into account.       Assuming that the utility functions at both stages are concave and continuously differentiable, we can derive the equivalent  two-stage stochastic equilibrium model  by the first-order optimality condition, which is a special case of the two-stage stochastic variational inequalities (SVI).

Variational inequalities play a central role in operations research and optimization, which model equilibrium problems in engineering and economics \cite{FePa97} and  represent the optimality conditions of the optimization problems \cite{RW-VA98}. Stochastic variational inequalities (also called  one-stage SVI ) \cite{ChFu05,ES11,GuOzRo99,GRS07,GF10,KaSKim11,KaSKim13,RaSha11} consider a
here-and-now solution in an uncertain environment using expectation of random functions. Compared with
the one-stage SVI, the two-stage SVI consider a pair of a here-and-now solution and a wait-and-see
solution, which  have inherently dynamic components that involve uncertain information and whose solutions depend on the outcome of random events to be revealed in the future.
 In the last few years, the two-stage SVI and multi-stage SVI have attracted considerable attention
\cite{ChPoWe17,ChSunXu18,ChShSun18,ChWe18,ChWeZh12,JSWC18,PaSeSh17,RoWe17,RoSun18a}.
 Chen, Pong and Wets \cite{ChPoWe17},  and Rockafellar and Wets \cite{RoWe17} introduced the
 the two-stage SVI and multi-stage SVI with examples including  the first-order optimality conditions for stochastic programs, Walras equilibrium problems with an incomplete financial market
and stochastic Wardrop traffic flow equilibrium problems. The existence and uniqueness of solutions of
the two-stage SVI, and convergence of its sample average approximation have been investigated in \cite{ChShSun18,ChSunXu18}. The Progressive Hedging Algorithm (PHA) \cite{RoWe91} for scenarios and policy aggregation in optimization under uncertainty has been extended to solve monotone two-stage SVI in \cite{RoSun18a}.

Our contributions in this paper are twofold. We develop new optimization theory and algorithms for the two-stage SVI arising from an oligopoly market. Furthermore, we apply the new optimization results to show the existence of a Nash-equilibrium for oil production game in an uncertain environment. In particular, we  model the game of crude oil production with a few agents as the two-stage SVI, and  show the existence and uniqueness of  solutions of  the two-stage SVI. Moreover, we develop a fast algorithm to find a solution of the two-stage SVI, which is used to analyze and forecast oil market share under the CODVI-19 pandemic.

This paper is organized as follows. In section 2, we present the two-stage SVI arising from two-stage stochastic games. In section 3, we show the existence and uniqueness of the solution of the two-stage SVI and provide perturbation error bounds for the
solution.  In section 4,  we develop a fast algorithm called Alternation Block Algorithm (ABA) to find a solution of
the two-stage stochastic equilibrium. The computational cost of the ABA is much less than the popular  PHA. We show outperformance of the ABA over the PHA using large-scale numerical examples.
In section 5, we apply the two-stage SVI, our new theory and algorithm to study the impact and the production responses of oil producers during COVID-19 on oil market share. Our numerical results with real market data from January 2019 to May 2020 show the efficiency of our methods for modelling oil market share.

 \section{Two-stage stochastic quadratic games}

We consider an oligopolistic market where $J$ agents compete to supply a
homogeneous product noncooperatively in the future. Each agent needs to make a decision on the quantity of production based on the anticipated future market demand and supply  and  other agents' decision in an uncertain environment. The uncertainties are represented by a random variable $\xi : \Omega \to \Xi\subset \mathbb{R}^m$ defined in the probability space $(\Omega, {\cal F}, {\cal P})$ with a support set $
\Xi$ and the space ${\cal Y}$ of measurable functions from $\Xi$ to $\mathbb{R}^J$.

We define the variables and functions for each agent $i$, $i=1,\ldots, J$.

$x_i\in \mathbb{R}$, the production quantity

$\theta_i: \mathbb{R}\to \mathbb{R}$, the cost function of the production

%$c_i\in \mathbb{R}$, the limit capacity to product quantity $x_i$

$y_i(\xi) \in \mathbb{R}$, the supply quantity

$\varphi_i:\mathbb{R}\times \Xi  \to \mathbb{R}$, the cost function for supplying a quantity

The uncertain market
supply and demand  are characterized by a random total supply $\eta(\xi)=\sum^J_{i=1}y_i(\xi)$ to the market and a random inverse
demand function $p(\eta(\xi),\xi): \mathbb{R}\times \Xi \to \mathbb{R}$.  Here
 $p(\eta(\xi),\xi)$ can be regarded as the spot price of trading in the future.

    Each agent aims to maximize its profit and make its decision in two stages.
      The first stage is to make optimal decision on production quantity based on the average of random demand, and the second stage is to make optimal decision on supply quantity based on
      the observation of the uncertainty in the future. The first stage decision is called ``here-and-now'' and the second stage decision is called ``wait-and-see''.

For each realization of random variable $\xi$ and a nonnegative vector $x \in \mathbb{R}^J$,   agent $i$  wants to find an optimal decision $y_i(\xi)$ by solving the following problem
 \begin{eqnarray}\label{eq2}
\begin{array}{cl}
\displaystyle F_i(x,\xi):=\max_{y_i(\xi)}  &  p(y_i(\xi)+y_{-i}(\xi),\xi)y_i(\xi) -\varphi_i(y_i(\xi), \xi)   \\
\hspace{0.9in} \mbox{s.t.}                    &0\leq y_i(\xi) \leq x_i,  \quad \,\, {\rm a.e.}\,\,  \xi \in \Xi.
\end{array}
\end{eqnarray}
Moreover, agent $i$ has to make an optimal decision before knowing the future events by solving
  \begin{eqnarray}\label{eq1}
\begin{array}{cl}
\displaystyle \max_{x_i}  & \bbe[F_i(x,\xi)] - \theta_i(x_i, x_{i-1})   \\
\mbox{s.t.}                    &0\leq x_i.
\end{array}
\label{eq:SNash}
\end{eqnarray}
Here  $x_{-i}$ and $y_{-i}$ are
decision variables of the agents other than the agent $i$.

Problem (\ref{eq1}) is the first stage, and problem (\ref{eq2}) is the second stage of the two-stage stochastic games \cite{ChShSun18,ChSunXu18,JSWC18,RaX11}.  We call $(x^*, y^*) \in \mathbb{R}^J \times {\cal Y}$ an optimal solution of the two-stage stochastic games, if for $i=1,\ldots, J$,
$(x_i^*, y^*_i(\xi))$ is the optimal solution of the two-stage optimization problem
  \begin{eqnarray}\label{eq11}
\begin{array}{cl}
\displaystyle \max_{x_i}  & \bbe[F_i(x_i, x_{-i}^*,\xi)] - \theta_i(x_i,x_{i-1}^*)   \\
\mbox{s.t.}                    &0\leq x_i,
\end{array}
\label{eq:SNash}
\end{eqnarray}
 \begin{eqnarray}\label{eq22}
\begin{array}{cl}
\displaystyle F_i(x_i, x_{-i}^*,\xi):=\max_{y_i(\xi)}  &  p(y_i(\xi)+y_{-i}^*(\xi),\xi)y_i(\xi) -\varphi_i(y_i(\xi), \xi)   \\
\hspace{1.1in} \mbox{s.t.}                    &0\leq y_i(\xi) \leq x_i,  \quad \,\, {\rm a.e.}\,\,  \xi \in \Xi.
\end{array}
\end{eqnarray}

%By the assumptions on the functions $p,\varphi_i$, problem (\ref{eq22}) is a strongly concave optimization problem and thus
%has a unique solution for any fixed $x$ and $\xi$. Moreover,
 %$F(x,\xi)$ is continuously differentiable w.r.t. $x_i$ for $x_i>0$ and
%\bgeq
%\nabla_{x_i}F_i(x,\xi) = \frac{L(y_i(\xi),\lambda_i(\xi),s_i(\xi),\x_i)}{dx_i}= s_i(\xi),
%\edeq
%where
%$$
%L(y_i(\xi),\lambda_i(\xi),s_i(\xi),\x_i) := p(y_i(\xi)+y_{-i}(\xi),\xi)y_i(\xi) -\varphi_i(y_i(\xi), \xi) +\lambda_i(\xi)y_i(\xi) - s_i(\xi)(y_i(\xi)- x_i).
%$$
%
%
%
%
%The Karush-Kuhn-Tucker (KKT) conditions  of problems (\ref{eq1})-(\ref{eq2}) derive the following two-stage stochastic  complementarity problem
%\bgeqn\label{eq:twogame}
%0\leq \begin{pmatrix}
%x_i\\
%y_i(\xi)\\
%s_i(\xi)
%\end{pmatrix} \perp
%\begin{pmatrix}
% \nabla_{x_i}\theta_i(x_i, x_{i-1})- \bbe[s_i(\xi)] \\
%g_i(y_i(\xi),y_{-i}(\xi),\xi)+ s_i(\xi)\\
%x_i - y_i(\xi)
%\end{pmatrix}\geq 0, \;\; &\inmat{for} \;\;\inmat{a.e.}\;\; \xi\in \Xi,  \;\; i=1,\ldots, J,
%\edeqn
%where $s_i$ is the Lagrange multiplier and
%$$
%g_i(y_i(\xi),y_{-i}(\xi),\xi) =
%-p(y_i(\xi)+y_{-i}(\xi),\xi)- y_i(\xi)p_q'(y_i(\xi)+y_{-i}(\xi),\xi)+\varphi_i'(y_i(\xi), \xi).
%$$

In this paper, we consider the case where the objective functions in (\ref{eq2})-(\ref{eq1}) are quadratic concave in the following forms
\begin{equation}\label{qp1}
\theta_i(x_i,x_{-i})=\frac{1}{2}c_ix_i^2 + a_ix_i+r_i x_i\sum^J_{j=1} x_j,
\end{equation}
\begin{equation}\label{qp2}
 \quad p(y_i(\xi)+y_{-i}(\xi),\xi)=\alpha(\xi) -\gamma(\xi) \sum^J_{j=1} y_j(\xi),
 \end{equation}
 \begin{equation}\label{qp3}
  \varphi_i(y_i(\xi), \xi)= \frac{1}{2}h_i(\xi) y_i^2(\xi)+\beta_i(\xi)y_i(\xi),
\end{equation}
where $c_i>0, \gamma(\xi)>0,  h_i(\xi)>0$, $a_i$, $r_i$, $\alpha(\xi)$ and $ \beta_i(\xi)$ are real given numbers.
In such setting, the function
$$
\begin{array}{cl}
\displaystyle F_i(x,\xi)= \max_{y_i(\xi)}& \displaystyle   ( \alpha(\xi) -\gamma(\xi) \sum^J_{j=1} y_j(\xi))y_i(\xi) -
\frac{1}{2}h_i(\xi) y_i^2(\xi)-\beta_i(\xi)y_i(\xi)\\
\hspace{0.7in} \mbox{s.t.}                    &0\leq y_i(\xi) \leq x_i,  \quad \,\, {\rm a.e.}\,\,  \xi \in \Xi
\end{array}$$
is  continuously differentiable with respect to (w.r.t.)  $x_i$ for $x_i>0$ and
$$\nabla_{x_i}F_i(x,\xi)
% = \frac{L(y_i(\xi),\lambda_i(\xi),s_i(\xi),x_i)}{dx_i}
= s_i(\xi),
$$
%where
%$$
%L(y_i(\xi),\lambda_i(\xi),s_i(\xi),x_i)$$
%$$=(\alpha(\xi) -\gamma(\xi) \sum^J_{j=1} y_j(\xi))y_i(\xi) -
%\frac{1}{2}h_i(\xi) y_i^2(\xi)-\beta_i(\xi)y_i(\xi) +\lambda_i(\xi)y_i(\xi) - s_i(\xi)(y_i(\xi)- x_i).
%$$
where
% $\lambda_i(\xi)\ge 0$ and
$s_i(\xi)\ge0 $ are Lagrange multipliers for the constraints
%$y_i(\xi)\ge 0$ and
$y_i(\xi) \leq x_i$.
If $x_i=0$, then the optimal solution of the optimization problem (\ref{eq22}) is $y_i(\xi)\equiv 0$
for any $\xi \in \Xi$, and the subdifferential $ \partial_{x_i} F_i(x,\xi)$ of $ F_i(x,\xi)$ is
$$\{s_i(\xi) \, |\,
s_i(\xi)\ge \max(0, \alpha(\xi) -\beta_i(\xi) -\gamma(\xi)e^Ty(\xi)), \, \bbe[s_i(\xi)]\le r_i e^Tx +a_i\}.$$

The Karush-Kuhn-Tucker (KKT) conditions  of problems (\ref{eq2})-(\ref{eq1}) with the functions defined in (\ref{qp1})-(\ref{qp3}) derive the following two-stage stochastic  linear complementarity problem (LCP)
\begin{equation}\label{eq:LCP}
0\leq \begin{pmatrix}
x\\
y(\xi)\\
s(\xi)
\end{pmatrix} \perp
\begin{pmatrix}
 (C+r e^T)x- \bbe[s(\xi)]+a \\
(H(\xi) +\gamma(\xi)ee^T)y(\xi)+ s(\xi) + \rho(\xi)\\
x - y(\xi)
\end{pmatrix}\geq 0,
\end{equation}
for almost every $\xi\in \Xi$,  where
$$x=(x_1,\ldots, x_J)^T, \quad a=(a_1,\ldots, a_J)^T, \quad  r=(r_1,\ldots, r_J)^T,$$
$$y(\xi)=(y_1(\xi), \ldots, y_J(\xi))^T, \quad
s(\xi)=(s_1(\xi), \ldots, s_J(\xi))^T,$$
$$
\rho(\xi)=(-\alpha(\xi)+\beta_1(\xi), \ldots, -\alpha(\xi)+\beta_J(\xi))^T,$$
$$ C= {\rm diag} (c_1+r_1,\ldots, c_J+r_J), \quad  H(\xi)={\rm diag}(h_1(\xi)+\gamma(\xi),\ldots, h_J(\xi)+\gamma(\xi)),$$
 and $e \in \mathbb{R}^J$ is the vector with all elements being 1.

 When the matrix $C+ re^T$ is positive definite, problems (\ref{eq2})-(\ref{eq1}) with the functions defined in (\ref{qp1})-(\ref{qp3}) are equivalent to problem (\ref{eq:LCP}) in the sense that
 if $(x^*, y^*(\cdot))$ is a solution of problems (\ref{eq2})-(\ref{eq1}), then there is $s^*(\cdot)$
 such that $(x^*, y^*(\cdot), s^*(\cdot))$ is a solution of (\ref{eq:LCP}); conversely, if
 $(x^*, y^*(\cdot), s^*(\cdot))$ is a solution of (\ref{eq:LCP}), then $(x^*, y^*(\cdot))$ is a solution of problems (\ref{eq2})-(\ref{eq1}). A sufficient condition for the matrix $C+ re^T$ being  positive definite is
\begin{equation}\label{cr}
c_i+ 2r_i > \frac{1}{2}\sum^J_{j\neq i}|r_j+r_i|,\quad
 i=1,\ldots, J.
 \end{equation}
  Condition (\ref{cr}) implies that $C+\frac{1}{2}(re^T+er^T)$ is a symmetric diagonally dominate matrix with positive diagonally elements and thus a positive definite matrix.

Let
 \begin{equation}\label{Mv}
  v(\xi)= \begin{pmatrix}
y(\xi)\\
s(\xi)
\end{pmatrix}\in \mathbb{R}^{2J}, \quad \,
M(\xi)=\begin{pmatrix}
H(\xi)+\gamma(\xi)ee^T  & \,\,  I\\
-I  & 0
\end{pmatrix} \in \mathbb{R}^{2J\times 2J},
%  \quad \,  B=\begin{pmatrix}
%0 & I
%\end{pmatrix} \in \mathbb{R}^{J\times 2J}.
\end{equation}

It is easy to see that for any  $x \in \mathbb{R}^J,
v \in \mathbb{R}^{2J},$
 \begin{equation}\label{PD}
(x^T, v^T)Q(\xi)\begin{pmatrix}
x \\
v
\end{pmatrix}=(x^T, v^T)\frac{1}{2}(Q(\xi)+Q(\xi)^T)\begin{pmatrix}
x \\
v
\end{pmatrix}\ge 0,
\end{equation}
where $ Q(\xi) = \begin{pmatrix}
C+re^T & -B \\
B & M(\xi)
\end{pmatrix} \in \mathbb{R}^{3J\times 3J}$  and $B=\begin{pmatrix}
0 & I
\end{pmatrix} \in \mathbb{R}^{J\times 2J}.$

Note that we cannot apply \cite[Proposition 2]{ChSunXu18} for the existence of solutions of the two-stage stochastic LCP (\ref{eq:LCP}), since the condition
 that there exists a positive continuous function $\kappa(\xi)$ with $\bbe[\kappa(\xi)]<\infty $ such that
\begin{equation}\label{PD}
(x^T, v^T)Q(\xi)\begin{pmatrix}
x \\
v
\end{pmatrix}\ge \kappa(\xi)(\|x\|^2+\|v\|^2),  \quad \forall x \in \mathbb{R}^J, \,
v \in \mathbb{R}^{2J},
\end{equation}
fails for any $x=0, v=(0, s)^T$, $s\in \mathbb{R}^J$ and $s\neq 0$.

%Moreover, we cannot apply \cite[Corollary 12.2]{HaKoSc05} to have the existence of solution of the two-stage stochastic LCP (\ref{eq:LCP}), since the function $G:\mathbb{R}^J\times {\cal V}\times {\cal V}\to \mathbb{R}^J\times {\cal V}\times {\cal V}$ defined by
%$$G(x,y,s)=\begin{pmatrix}
% (C+r e^T)x- \bbe[s(\xi)]+a \\
%(H(\xi) +\gamma(\xi)ee^T)y(\xi)+ s(\xi) + \rho(\xi)\\
%x - y(\xi)
%\end{pmatrix}$$
%is not weakly coercive.
%\footnote{$\mathcal{G}$ is weakly coercive if
%there exists $(x_0, y_0,s_0)\in \mathbb{R}^J\times \times {\cal V}\times {\cal V}$ such that
%$$
%\left\langle G(x,y,s),  \begin{pmatrix}
%x\\
%y\\
%s
%\end{pmatrix}\right\rangle= x^T (C+r e^T)x
%+ \bbe[y(\xi)^T(H(\xi) +\gamma(\xi)ee^T)y(\xi) +s(\xi)]+\bbe[s(\xi)^T(x - y(\xi))]
%\to \infty \;\;\inmat{as}\;\; \|x\| + \|y \|+\|s\| \to\infty \;\;\inmat{and}\; (x, y,s)\in \mathbb{R}^J\times \times {\cal V}\times {\cal V}.
%$$
%}

In the next section, we will establish the existence of solutions of the two-stage stochastic LCP (\ref{eq:LCP}) with a finite number of realizations of the random variable.

%The existence of solutions of the two-stage stochastic games (\ref{eq1})-(\ref{eq2}) withe the quadratic functions in (\ref{qp1})-(\ref{qp2})
%%and the two-stage stochastic LCP (\ref{eq:LCP})
% is guaranteed by the strongly concavity of these functions. Moreover,
%
%{\bf Remark} In the rest of this paper, we will assume the component $x^*$ of solution is positive.
%If $x^*_i=0$, then by the constraint of the second stage, we have $y_i^*(\xi)=0$ for a.e. $\xi \in \Xi$. In this case, we can remove agent $i$ from
%the game.

% Hermitian strictly diagonally dominant matrix with real positive diagonal entries is positive definite

\section{Existence, uniqueness and robustness of solutions }

%In this section, we will study the link of parameters $C$, $a$, $r$, $H(\xi)$, $\gamma(\xi)$ and $\rho(\xi)$ to the solution $(x^*,y^*(\xi)$) of problem (\ref{eq:LCP}).

   We first  consider the following LCP$(M(\xi),q(x,\xi))$
 $$
0\leq v(\xi) \perp M(\xi) v(\xi)  + q(x,\xi)
\geq 0,
$$
for a fixed $x\ge0$ and $\xi \in \Xi$, where $M(\xi)$ and $v(\xi)$ are defined in (\ref{Mv}) and $ q(x,\xi)= \begin{pmatrix}
\rho(\xi)\\
x
\end{pmatrix}\in \mathbb{R}^{2J}.$

The matrix $M(\xi)$ is positive semi-definite, but not positive definite, since
$$v(\xi)^TM(\xi)v(\xi)=\frac{1}{2}v(\xi)^T(M(\xi)+M(\xi)^T)v(\xi)=y(\xi)^T(H(\xi)+\gamma(\xi)ee^T)y(\xi)\ge 0$$
and $ v(\xi)^TM(\xi)v(\xi)=0$ if $y(\xi)=0$ and $s(\xi)\neq0$.

  Note that the matrix $M(\xi)$ is not a $P$-matrix. Thus the existence and uniqueness of the solution of the LCP$(M(\xi),q(x,\xi))$ and the error bound in \cite{Chen_Xiang1,Chen_Xiang2,CoPaSt92} cannot be guaranteed for any $q(x,\xi)$.
For example, if there is $i$ such that $x_i< 0$, then the LCP$(M(\xi),q(x,\xi))$ does not have a solution; if  $x=0$, then any $v(\xi)$ with $y(\xi)=0$ and $s(\xi)\ge -\rho(\xi)$ is a solution of
the LCP$(M(\xi),q(x,\xi))$.
%On the other hand, the matrix $M(\xi)$ is nonsingular.

Let SOL$(M(\xi),q(x,\xi))$ be the solution set of the LCP$(M(\xi),q(x,\xi))$. Since $M(\xi)$ is a
positive semi-definite matrix, SOL$(M(\xi),q(x,\xi))$ is a convex set \cite{CoPaSt92}.
 A vector
$v^*(\xi)$  is called a least norm solution of  the LCP$(M(\xi),q(x,\xi))$ if
it is the solution of the optimization problem
$$ \min \|v\|^2  \quad {\rm s.t.} \, v\in {\rm SOL}(M(\xi),q(x,\xi)),$$
where $\|\cdot\|$ is the Euclidean norm.

The following lemma gives the form of the least norm solution of the LCP$(M(\xi),q(x,\xi))$ and shows that it is the unique solution  of the LCP$(M(\xi),q(x,\xi))$ when $x>0$.
\begin{lemma}
For any $x\ge 0$ and $\xi \in \Xi$, the LCP$(M(\xi),q(x,\xi))$ has the least norm solution $v^*(\xi)=(y^*(\xi), s^*(\xi))^T$  with
\begin{equation}\label{fixed}
y^*(\xi)=\Pi_{[0, x]}(y^*(\xi)-\rho(\xi)-  (H(\xi)+\gamma(\xi)ee^T)y^*(\xi)),
\end{equation}
\begin{equation}\label{sols}
s^*(\xi)=\max(0, -\rho(\xi)- (H(\xi)+\gamma(\xi) ee^T)y^*(\xi)),
\end{equation}
where $\Pi_{[0, x]}(z)$ is the projection from $z$ to the set $[0,x]$.
Moreover, the least norm solution is the unique solution of the LCP$(M(\xi),q(x,\xi))$ if $x>0$.
\end{lemma}

Lemma 3.1 shows that for a fixed vector $x\ge 0$, the least norm  solution of the LCP$(M(\xi),q(x,\xi))$ is
uniquely defined by $H(\xi)$, $\rho(\xi)$ and $\gamma(\xi)$. In real applications, noise may be presented in the data set and  we  consider the perturbation bound of
the solution regarding noise in $H(\xi)$, $\rho(\xi)$ and $\gamma(\xi)$.

 Let us consider the following
 %quadratic program
%\begin{equation}\label{qp3}
%\begin{array}{cc}
%\min & \frac{1}{2} u^T(\bar{H} + \bar{\gamma}ee^T)u +\bar{\rho}^Tu\\
%{\rm s.t.} &  0\le u\le x,
%\end{array}
%\end{equation}
%and the corresponding
LCP$(\bar{M}, \bar{q})$ with
$$ \bar{M}=\begin{pmatrix}
\bar{H}+\bar{\gamma}ee^T  & \,\,  I\\
-I  & 0
\end{pmatrix},  \quad \,  \bar{q}=\begin{pmatrix}
\bar{\rho}\\
x
\end{pmatrix}.
$$

From Lemma 3.1,  the LCP$(\bar{M}, \bar{q})$ with $x\ge 0$ has the least norm solution $\bar{v}=(\bar{u},\bar{t})^T$ in the following form
$$
\bar{u}=\Pi_{[0, x]}(\bar{u}- \bar{\rho} -  (\bar{H}+\bar{\gamma}  ee^T)\bar{u})
\quad {\rm and} \quad
\bar{t}=\max(0,
-\bar{\rho} - (\bar{H}+\bar{\gamma} ee^T)\bar{u}).
$$

The following theorem provides the distance $\|y^*(\xi)-\bar{u}\|$  regarding the noise in data  $\rho(\xi), H(\xi), \gamma(\xi).$
\begin{theorem} Suppose that $x\ge 0$ and let $\Gamma>0$ such that $\Gamma \gamma(\xi)\ge 1$. Then the least norm solution of the LCP$(M(\xi), q(x,\xi))$ is continuous with respect to $H(\xi)$, $\rho(\xi)$ and $\gamma(\xi)$. Moreover, we have the following perturbation error bound
\begin{equation}\label{error_bound}
\|y^*(\xi)-\bar{u}\|\le \Gamma(\|\rho(\xi)-\bar{\rho}\|+\|x\|\|H(\xi)-\bar{H}\|+J\|x\|\|\gamma(\xi)-\bar{\gamma}\|).
\end{equation}
\end{theorem}

\section{A new Alternating Block Algorithm (ABA)}

In this section, we consider how to efficiently solve   problem (\ref{eq:LCP}) with a finite number of realizations of the random variable. For $\Xi^\nu=\{\xi_1,\ldots, \xi_\nu\}$, the probability for each realization is $1/\nu$.
For simplicity, we set $B=(0, I)\in \mathbb{R}^{J\times 2J}$, $\varrho(\xi)=(\rho(\xi), 0)^T
\in \mathbb{R}^{2J}$, $n=J+2J\nu$ and
$$
%\mathbf{v}=\left(\begin{array}{c}
%x\\
%v(\xi_1)\\
%\vdots\\
%v(\xi_\nu)
%\end{array}
%\right) \in \mathbb{R}^n, \quad
\mathbf{M}=\left(\begin{array}{cccc}
C+re^T &   -\frac{1}{\nu}B & \dots &  -\frac{1}{\nu}B\\
B^T & M(\xi_1) &  & \\
\vdots& & \ddots &\\
B^T    &     &      & M(\xi_\nu)
\end{array}
\right) \in \mathbb{R}^{n\times n  }, \quad \mathbf{q}=\left(\begin{array}{c}
a\\
\varrho(\xi_1)\\
\vdots\\
\varrho(\xi_\nu)\\
\end{array}
\right)\in \mathbb{R}^{n}.$$

In such setting, problem (\ref{eq:LCP}) is the standard LCP$(\mathbf{M}$,$\mathbf{q})$ and
the  progressive hedging algorithm  for the LCP$(\mathbf{M}$,$\mathbf{q})$ is as
follows.

  \parbox{\textwidth}{%
{\bf  Algorithm 1:  Progressive Hedging Algorithm (PHA)} \cite{RoSun18a}

{\bf Step 0.} Given an initial point $x^0 \in \mathbb{R}^J$, let $x_\ell^0=x^0\in \mathbb{R}^J, v_\ell^0 \in \mathbb{R}^{2J}$ and $w_\ell^0 \in \mathbb{R}^J$, for $\ell=1,\ldots,\nu$,   such that $\frac{1}{\nu}\Sigma_{\ell=1}^\nu w_\ell^0=0.$ Set the initial point $z^0=(x^0, v_1^0,\ldots,v_\nu^0)^T$. Choose a step size $t>0.$  Set $k=0.$

%{\bf Step 1.} If the point $z^k$  satisfies the  condition
%$$ \|\min(z^k, H^\epsilon z^k+\bar{q}) \|\leq 10^{-6}, $$
%output the solution $z^k$ and terminate the algorithm; otherwise,  go to Step 2.

{\bf Step 1.}  For $\ell=1,\ldots,\nu$, find $(\hat{x}^k_\ell, \hat{v}_\ell^k)$ that solves the LCP
\begin{equation}\label{PHM-subproblem}
  \begin{array}{l}
  0\leq x_\ell \bot (C+re^T) x_\ell -B v_\ell+a +w_\ell^{k} +t(x_\ell-x_\ell^k)  \geq 0,  \\
  0\leq v_\ell \bot  B^T x_\ell+ M(\xi_\ell) v_\ell  + \varrho(\xi_\ell)+t(v_\ell-v_\ell^k)  \geq 0.
\end{array}
\end{equation}
 Let $\bar{x}^{k+1}=\frac{1}{\nu} \sum_{\ell=1}^\nu \hat{x}_\ell^k$,  and for $\ell=1,\ldots,\nu,$ update
\begin{equation}
  x_\ell^{k+1}=\bar{x}^{k+1}, \,\, v_\ell^{k+1}=\hat{v}_\ell^{k}, \,\, w_\ell^{k+1}=w_\ell^{k}+t(\hat{x}_\ell^{k}-x_\ell^{k+1}),\notag
\end{equation}
to get point $z^{k+1}=(\bar{x}^{k+1},v_1^{k+1},\ldots,v_\nu^{k+1})^T$.

{\bf Step 2.} Set $k:=k+1;$  go back to Step 1.\vspace{2mm}
  }%

\begin{theorem}
The LCP$(\mathbf{M}$,$\mathbf{q})$ has at least one solution $\mathbf{v}^*$, and has at most one solution with $x^*>0$. Moreover, the sequence $\{z^k\}$ generated by the PHA converges to a solution of the LCP$(\mathbf{M}$,$\mathbf{q})$.
\end{theorem}

The PHA has been widely  used for solving stochastic optimization problems and stochastic variational inequalities. Due to the special structure of the $\mathbf{M}$, we propose a new algorithm for solving the  LCP$( \mathbf{M}$,$\mathbf{q})$, which is called Alternating Block Algorithm (ABA).
At each iteration of the PHA, we need to solve $\nu$ linear complementarity problems in $3J$ dimension.
 At each iteration of the ABA, we only need to solve $\nu$ strongly convex quadratic programs with a simple constraint in $J$ dimension and one linear complementarity problem in $J$ dimension. Hence the computation cost of the ABA is much less than that of the PHA.

  \parbox{\textwidth}{%
{\bf  Algorithm 2:  Alternating Block Algorithm (ABA)}

{\bf Step 0.} Given an initial point $x^0 \in \mathbb{R}^J$ with $x^0\ge 0$.  Set $k=0.$

%{\bf Step 1.} If the point $z^k$  satisfies the  condition
%$$ \|\min(z^k, H^\epsilon z^k+\bar{q}) \|\leq 10^{-6}, $$
%output the solution $z^k$ and terminate the algorithm; otherwise,  go to Step 2.

{\bf Step 1.}  For $\ell=1,\ldots,\nu$, find $y^k_\ell$ that solves quadratic program
\begin{equation}\label{AB-subproblem}
 \begin{array}{cc}
\min & \frac{1}{2} y^T(H(\xi_\ell) + \gamma(\xi_\ell)ee^T)y +\rho(\xi_\ell)^Ty\\
{\rm s.t.} &  0\le y \le x^k.
\end{array}
\end{equation}
 Set
% \begin{equation}\label{AB-sols}
%(s^k_\ell)_i=\left\{\begin{array}{ll}
%0 & \quad {\rm if} \, (y^k_\ell)_i< x_i^k \\
% -(\rho(\xi_\ell)+ (H(\xi_\ell)+\gamma(\xi_\ell) ee^T)y^k_\ell)_i & \quad {\rm if } \, (y^k_\ell)_i= x_i^k.
%\end{array}
%\right.
%\end{equation}
\begin{equation}\label{AB-sols}
s^k_\ell=\max (0,  -\rho(\xi_\ell)- (H(\xi_\ell)+\gamma(\xi_\ell) ee^T)y^k_\ell).
\end{equation}
{\bf Step 2.} Find $x^{k+1}$  that solves the LCP problem
$$0 \le  x \bot  (C+re^T)x -\frac{1}{\nu} \sum^\nu_{\ell=1} s^k_\ell + a \ge 0.$$
Set $k:=k+1;$  go back to Step 1.\vspace{2mm}
  }%

Let $C_2\in \mathbb{R}^{J\times J}$ be the symmetric positive definite matrix  such that
$C_2^TC_2=C+\frac{1}{2}(re^T+er^T).$ Let
$$\sigma_\ell= \max_{K\subseteq\{1,\ldots, 2J\}} \|M(\xi_\ell)_{K}^{-1}\|, \quad \sigma=\frac{1}{\nu}\sum^\nu_{\ell=1} \sigma_\ell,$$
where $M(\xi_\ell)_{K}$ is the nonsingular principal submatrix of $M(\xi_\ell)$ whose entries of $M(\xi_\ell)$ are indexed by the set $K\subseteq  \{1,\ldots, 2J\}.$

\begin{theorem}
Suppose $\|C_2^{-1}\|^2\sigma <1$ and  the first $J$ component $x^*$ of the solution of the LCP$(\mathbf{M}$,$\mathbf{q})$ is positive.  Then there is neighborhood of $x^*$ such that for any initial point $x^0$ in the neighborhood, the sequence generated by the ABA converges to the solution of the LCP$(\mathbf{M}$,$\mathbf{q})$.
\end{theorem}

Now we use randomly generated problems to compare the performance of  the PHA and ABA.
With uniform distribution, we randomly generate each element of vectors $a, c,  \bar{\beta} $
on $[0, 1]$, $\bar{h}$ on $[2,3]$, and numbers $\bar{\gamma}$ on $[0,0.5]$ and $ \bar{\alpha}$ on $[5,10]$. Let $r=0.5e$ and
$C_{ii}=10+c_{i} + (r^Te + (J-2)r_i).$
We generate sample $\xi_1,\xi_2,\ldots, \xi_\nu$ from uniform distribution on $[1,2]$ and set
\begin{align*}
\beta(\xi_\ell)=\xi_\ell \bar{\beta},\, h(\xi_\ell)=\xi_\ell \bar{h},\,\gamma(\xi_\ell)=\xi_\ell\bar{\gamma}, \,
\alpha(\xi_\ell)=\xi_\ell \bar{\alpha}, \quad \ell=1,\ldots, \nu.
\end{align*}
Then, we let
\begin{align*}
\rho(\xi_\ell)=(-\alpha(\xi_\ell)e+\beta(\xi_\ell)),\,\,    \,\, H(\xi_\ell)={\rm diag}(h(\xi_\ell)+\gamma(\xi_\ell)e), \\
M(\xi_\ell)=
 \left(
   \begin{array}{cc}
     H(\xi) +\gamma(\xi_\ell)ee^T & I \\
     -I & 0 \\
   \end{array}
 \right), \,\,\,\, \varrho(\xi_\ell)=\left(
   \begin{array}{c}\rho(\xi_\ell)\\
   0
   \end{array}\right),\quad \ell=1,\ldots, \nu.
\end{align*}
%We take $\epsilon=10^{-12}$ in order to compute a least-norm solution of ${\rm LCP}(\mathbf{M}, \mathbf{q})$.

We terminate the two algorithms when one of the following three stop criteria is met:
The number of iterations reaches  400,
or $\|\mathbf{v}_k-\mathbf{v}_{k-1}\|\leq 10^{-6}$, or
%\begin{equation}\label{res-cri}
$${\rm Res}:=\|\min(\mathbf{M}\mathbf{v}+\mathbf{q},\mathbf{v})\|\leq 10^{-6}.
$$
%\end{equation}

%If \eqref{res-cri} is satisfied, we say that $\mathbf{v}$ is a solution of $LCP(\mathbf{M},\mathbf{q})$.
 The initial points for PHA and ABA are $\mathbf{v}_0=((-(C+re^T)^{-1}a)_{+},0,\ldots,0)^T$ and
  $x^0=-(C+re^T)^{-1}a)_{+}$, respectively. The step size $t$ in the PHA is set to $t=1$.
%$$\|\min ( C+re^T) x^k -\frac{1}{\nu}\sum^\nu_{\ell=1} s_\ell^k)\|^2
%+ \frac{1}{\nu}\sum^\nu_{\ell=1} \|\min (M_\ell v_\ell^k + B^Tx^k +q_\ell, v_\ell^k) \|^2 \le 10^{-10}$$
%\begin{center}

We choose $J=5, 15, 20$,  and increase the sample size $\nu$ from 5 to 1000. The dimension of the corresponding ${\rm LCP}(\mathbf{M},\mathbf{q})$ ranges
from 55 to 30015. For each  $J, \nu$, we randomly generate 10 problems following the  description above. The ABA and PHA methods are used to solve these 10 problems. The results reported in
 Table 1 are the average of iterations, cpu times(seconds), residuals and initial residuals.
 Table 1 shows that the ABA can solve the ${\rm LCP}(\mathbf{M},\mathbf{q})$ with all $J$ and $\nu$ efficiently.
Another advantage of the ABA is that the iteration numbers remain almost the same when $J$ and $\nu$ increase, and the cpu time
is roughly linearly increasing  as $\nu$ increases. Moreover Table 1 shows that
although the PHA can solve the problems with small $J$ and $\nu$, it fails to solve the problems with
large $J$ or $\nu$ within 400 iterations.
%For the problem with dimension 30015, PHA requires 100 seconds to perform 400 iterations reaching an accuracy of $10^{-3}$ only, while ABA do 11 iterations in 22 seconds to
%compute a solution with accuracy of $10^{-6}$.

\begin{table}\label{table-pha-aba}
\begin{tabular}{|c|c|c|c|c|c|c|c|}
  \hline
  \multicolumn{1}{|c|}{  }&\multicolumn{3}{|c|}{ABA}&\multicolumn{3}{|c|}{ PHA} &\multicolumn{1}{|c|}{  }\\
\hline
$J, \nu, J(2\nu+1)$ & iter &   CPU & Res&    iter &    CPU &Res & Initial Res\\
\hline
       5           5          55 &      15.60 &       0.02   &   6.40e-7 &    294.10 &       0.09   &   9.76e-7 &   5.37e+1 \\
     5          50         505 &      18.40 &       0.06   &   7.62e-7 &    342.90 &       0.58   &   9.91e-7 &   1.69e+2 \\
     5         100        1005 &      21.70 &       0.11   &   6.13e-7 &    371.20 &       1.18   &   1.23e-6 &   2.68e+2 \\
     5         500        5005 &      22.30 &       0.55   &   6.82e-7 &    379.00 &       5.92   &   2.25e-6 &   5.80e+2 \\
     5        1000       10005 &      22.00 &       1.11   &   6.37e-7 &    387.20 &      12.06   &   2.16e-6 &   8.38e+2 \\
\hline
    10           5         110 &      20.10 &       0.05   &   8.03e-7 &    383.10 &       0.22   &   5.36e-6 &   7.29e+1 \\
    10          50        1010 &      20.60 &       0.21   &   7.18e-7 &    399.20 &       1.66   &   1.64e-5 &   2.42e+2 \\
    10         100        2010 &      25.20 &       0.67   &   7.12e-7 &    400.00 &       3.63   &   4.22e-5 &   3.18e+2 \\
    10         500       10010 &      25.10 &       2.34   &   5.88e-7 &    400.00 &      17.75   &   9.66e-5 &   8.13e+2 \\
    10        1000       20010 &      23.00 &       6.31   &   5.95e-7 &    400.00 &      34.80   &   9.27e-5 &   1.01e+3 \\
\hline
    15           5         165 &      14.50 &       0.04   &   6.85e-7 &    400.00 &       0.40   &   1.21e-5 &   9.31e+1 \\
    15          50        1515 &      20.70 &       0.69   &   6.76e-7 &    400.00 &       3.66   &   1.74e-4 &   2.69e+2 \\
    15         100        3015 &      20.10 &       1.26   &   5.69e-7 &    400.00 &       9.97   &   1.62e-4 &   3.78e+2 \\
    15         500       15015 &      17.80 &       3.91   &   5.56e-7 &    400.00 &      48.55   &   3.11e-4 &   9.02e+2 \\
    15        1000       30015 &      21.60 &      11.31   &   1.02e-6 &    400.00 &     100.23   &   7.26e-4 &   1.41e+3 \\
\hline
\end{tabular}
\caption{Comparison of the PHA and the ABA \label{blackwhite}}
\end{table}
%\end{center}
%In the next section, we will use Algorithm 1 for our numerical test.

%\begin{theorem}
%The two-stage stochastic LCP (\ref{eq:LCP}) has a solution in $\mathbb{R}^J \times {\cal V}\times {\cal V}$. Moreover, $(x^*,y^*) \in \mathbb{R}^J \times {\cal V}$ is an optimal solution of the two-stage stochastic games if and only if
%$(x^*,y^*,s^*)\in \mathbb{R}^J \times {\cal V}\times {\cal V}$ is a solution of the two-stage stochastic LCP.
%\end{theorem}
%
%{\bf To be proved.}

\section{Impact of COVID-19 on oil market share}

In this section, we use the real data of oil price, demand and the market share of 14 major oil producers in the last 17 months to demonstrate the predicability of the two-stage stochastic LCP model with adaptive  parameters in the cost functions. Moreover, the model and the simulation results rationalize decisions made by major producers during COVID-19 pandemic.

Since its first identification, COVID-19 has become a pandemic across all continents of the globe.
To fight back the disease, large portions of social and commercial activities have been restricted if not suspended entirely.
It seemed, during the worst of the outbreak, the whole world had come to a halt, so did the consumption of oil.

This far reaching global event contributed to an already fragile market of oil.
%US joined the ranks of oil exporters in 2014 due to large scale extraction activities of shale oil, and the traditional oil producers responded passively by continuing production with no cuts. As a result, the crude oil price drop from over US\$110 in 2014 to below US\$ 30 in 2016 measured in Brent benchmark price. In September 2016, Saudi Arabia and Russia started to cooperate in managing the price of oil by cutting their production, creating an informal alliance of OPEC and non-OPEC traditional producers, commonly known as OPEC+. The agreement of production cut continued and by January 2020, OPEC+ had cut oil production by 2.1 million barrels per day comparing to its peak production.
As a result of the pandemic, demand in transportation and factory output fell, which leads to decrease in overall demand for oil.
Consequently, it caused oil prices to fall deeply. The last straw happened on 6 March 2020, when Russia rejected the demand of OPEC on further production cuts in response to the demand shrinkage. The oil price fell following the Russian announcement.

%It may be strange at first glance, since a further cut in production should, in theory, boost the price. Wouldn't that be beneficial to producers of oil? Not necessarily.
%High price seems to be attractive in short term, it may hurt in the long run.
%One typical example of such was triggered by the oil embargo from 1973-1974, although supply restriction leaded to price spikes in turn caused considerable economic damage for developed countries.
%When a similar situation presented itself during 1979-1980 oil crisis, strategy of supply restriction back-fired, and resulted in free-fall of the collective market share of all members of OPEC, from over 50{} in 1979 to below 30{} by the end of 1985.

\subsection{Modelling oil market share as a two-stage stochastic game}

We treat the oligopoly market of oil as a two-stage stochastic game where oil producers compete for profit by deciding the optimal production at the first stage of the game.
In particular, we consider 15 oil producers in the following list as 15 agents in the game.\\

1 Saudi Arabia,\, 2 Russia,\,  3 USA, \,
4 Iraq,\,
5 China,\,
6 Canada,\,
7  United Arab Emirates (UAE),\,  8 Iran,\,
9 Kuwait, \, 10	Nigeria, \, 11 Mexico, \, 12 UK, \,  13  Venezuela,	\,
14 Indonesia, \, 15 other.	 \\

Producer $i$ makes a decision on its oil production quantity $x_i$, based on its predicted further characteristics of the oil market at a later time, where the trading actually occurs. We suppose the trading occurs at the second stage where producer $i$ supplies part if not all of its produced quantity in the first stage to generate revenue. The spot price of the trading is uncertain at the time of production decision, and is mainly depended on the demand supply relation at the second stage. To be more precise, supply quantity $y_i(\xi)$ of producer $i$ is uncertain being further event at the time of production decision.

To simulate the effect of excessive supply on price, we adapt a simple supply and demand relation and express the inverse demand function as:
$$p(y(\xi),\xi) = p_0(\xi) -\gamma(\xi)\Big(\sum^J_{j=1} y_j(\xi) -D(\xi)\Big),$$
where $p_0(\xi)$ is the benchmark price and $\gamma(\xi)>0$ represents the negative effect on price if there is excessive supply of oil with respect to the observed demand $D(\xi)$ at the second stage.

Traditionally, for the production of every barrel of crude oil, producers need to explore oil fields before building the extraction site, even after the oil is extracted refinement and shipment require both time and labour not to mention the cost that involves. The nature of the oil production has been evolving with the technological advance which enables, e.g., extraction from oil sands.
There exist fundamental differences in energy infrastructures between traditional producers and oil sands producers, e.g., US, Canada \cite{Lazzaroni},
and it is expected to be reflected at the first stage in choosing production strategies.

%There are costs involved in both the production stage and the supply stage.  For the supply stage, the produced oil needs to be refined and transported to the destinations. Since the oil is produced a prior to the trading time the storage cost can also make an impact of adjusted price of each agent. In fact, one major reason that caused the WTI price to fall below US\$ zero per barrel was the lack of storage facilities, not only at Cushing where the oil is physically transferred but also over the entire global where hundreds of fully loaded oil tankers being stationary at sea with no destinations to go and oil pipes are taken as temporal storage. The composition of production cost is also complicated since it involves exploration, building facilities, extraction and many more steps.
In this paper, we assume that the costs in both stages are quadratic as expressed in (\ref{qp1}) and (\ref{qp2}) where the parameters $c_i, a_i, h_i(\xi), \gamma(\xi)$ are to be learned from market data.
%The parameter $\rho_i(\xi)$ can be considered as the adjusted price for each producer with the expression $\rho_i(\xi) = p_0(\xi)+\gamma(\xi)D(\xi)+w_i(\xi)$, where $w_i(\xi)$ represents the parameter of the linear part of second stage cost.
We also include an extra term $(r_i\sum^J_{j=1} x_j)x_i$ in the first stage to represent the strategic concern of producer $i$ in response to global production quantity
$\sum^J_{j=1} x_j$. Note that we have no restriction on the sign of $r_i$ and when it is positive it represents the fact that producer $i$ is willing to decrease its
production when the global production is high. Typical of such agents are often price setters, since the action would give a boost to the oil price and may be more profitable despite the production cuts.
In reality, Russia's refusal to production cuts agreement triggered huge volatility of the market during the COVID-19 pandemic.
To express it with our model, it means that Russia adapted its value of $r_i$ to be negative given the market prediction back in March 2020. That is to say, Russia made a ``squeeze'' strategy and decided that if the price drop caused by over supplying is under control, the maintained market share is potentially more important for long term profitability.
Other notable decisions of producers during the COVID-19 pandemic include that  Russia refused to production cut on 6th March 2020; Saudi Arabia offered price discount on 8th  March and
increased production; U.S. demanded production cuts on 2nd  April 2020;
OPEC and Russia cut production on 9th  April 2020.
%\begin{itemize}
%\item Russia's refusal to production cut on 6th of March 2020;
%\item Saudi Arabia's decision to offer price discount on 8th of March;
%\item Saudi Arabia's decision to increase production and Russia's decision of increase that followed;
%\item U.S. demanded production cuts 2nd of April 2020;
%\item OPEC and Russia's decision to cut production 9th of April 2020.
%\end{itemize}

We will use our model to show that most decisions of producers were reasonable from the prospective of the producers, since they had different tolerance on low oil prices, and they all attempted the best action for their own benefits.\\

\subsection{Numerical simulation  and parameter setting}
%Below is the real historical data of crude oil from the reliable sources.
The market data used in our study are obtained from the following sources.
\begin{description}
\item(i) Statistical Review of World Energy\footnote{https://www.bp.com/en/global/corporate/energy-economics/statistical-review-of-world-energy.html}, %  yearly data from 1984,
    latest publish in June 2020 by bp. Inc;

This data set provides average daily production quantities of major oil-producing countries, and the percentage of their production of the total production is regarded as their market shares respectively.

\item(ii) Oil Price Dynamics Report\footnote{https://www.newyorkfed.org/research/policy/oil\_price\_dynamics\_report}, weekly by
Federal Reserve Bank of New York;

This data set reports the crude oil price change and the supply/demand relation contributing to the change. It shows supply contribution, demand contribution and other contribution called residual contribution.

\item(iii) U.S. Energy Information Administration\footnote{https://www.eia.gov}, weekly and daily spot price of Brent.

\end{description}

For the model (\ref{eq2})-(\ref{eq1}) with (\ref{qp1})-(\ref{qp3}), parameters in the first stage are taken as follows:
%There are following model parameters to be determined using the above real data.
\begin{itemize}
\item $c_i$: This parameter represents the quadratic production cost of producer $i$. For the traditional producers, this contributes to the cost of exploration, site building and equipment setting up, etc. For oil sand producers, the financial cost can also be regarded as non-linear with respect to production quantity. However, all major producers should have similar scale of values to maintain profitable.
In our simulation of in-sample and out-of-sample,  we took
$$c_1=0.11/\Lambda_1, \,\, c_2=0.115/\Lambda_2, \,\, c_3=0.095/\Lambda_3,\,\, c_i=0.1/\Lambda_i, \,\, i=4,\ldots, 15$$
$$c_1=0.11/\Lambda'_1, \,\, c_2=0.115/\Lambda'_2, \,\, c_3=0.095/\Lambda'_3,\,\, c_i=0.1/\Lambda'_i, \,\, i=4,\ldots, 15$$
respectively, where $\Lambda_i$ and $\Lambda'_i$ are  market share of producer $i$ in the current month and previous month for simulation of  2019,  and market shares of January 2020 and December 2019 for simulation of  2020.
The data of market share are given in Tables 3-4.
 %so that producers with similar share are kept close in its $c_i$ values.

\item $a_i$: This parameter represents the linear production cost of producer $i$. It is widely agreed that traditional producers have very low unit cost of oil production. As for the oil sand producers, the unit cost is much higher.
In our simulation,
$$a_i=c_i, \, i=1,2,4,5,7,\ldots, 15, \,\, a_3=6c_3, \, a_6=2c_6.$$
 %is set to be equal $c_i$ for traditional producers while taking values of $10c_i$ for oil sand producers.

\item $r_i$: This parameter represents producer $i$'s  response to the total of production by all producers.    Before the pandemic, it was widely agreed that supply should be kept in accordance to the demand but little preference was taken till Russia's refusal to further production cuts. In our simulation,  $r_i=0$ for all producers in 2019, and $r_i$ were given in Table 2
    for different producers in 2020.
\end{itemize}

\begin{table}
\begin{center}
%\begin{tiny}
\begin{tabular}{|c|c|c|c|c|c|}
  \hline
2020	& Jan &	Feb &	Mar	& Apr &	May	\\
\hline
Saudi Arabia&	0{}	&	0.01{}	&	0{}	&	-0.022{}	&	0{}	\\
Russia	&-0.01{}	&	0{}	&	0{}	&	-0.008{}	&	0.01{}	\\
USA	&-0.01{}	&	0{}	&	-0.02{}	&	-0.04{}	&	0{}	\\
Iraq	&0{}	&	0{}	&	0{}	&	-0.01{}	&	0{}	\\
China&	0{}	&	0{}	&	0.01{}	&	-0.01{}	&	0{}	\\
Canada	&0{}	&	0{}	&	-0.02{}	&	-0.03{}	&	0{}	\\
UAE&	0{}	&	0{}	&	-0.02{}	&	-0.05{}	&	0{}	\\
Iran	&0{}	&	0{}	&	-0.02{}	&	-0.045{}	&	0{}	\\
Kuwait	&0{}	&	0{}	&	-0.01{}	&	-0.045{}	&	0{}	\\
Nigeria&	0{}	&	0{}	&	-0.01{}	&	-0.08{}	&	-0.06{}	\\
Mexico	&0{}	&	0{}	&	0{}	&	-0.065{}	&	-0.045{}	\\
UK	& 0{}	&	0{}	&	-0.1{}	&	-0.16{}	&	-0.08{}	\\
Venezuela	&0{} &0{}	&	-0.1{}	&	-0.23{}	&		-0.13{}	\\
Indonesia	&0{} &0{}	&	-0.1{}	&	-0.23{}	&		-0.13{}	\\
other	&-0.01{} &0{}	&	0.005{}	&	0.005{}	&		-0.005{}	\\
  \hline
\end{tabular}
%\end{tiny}
 \caption{Values of $r$ for numerical  simulation of January-May 2020}\label{table:share-2020}
 \end{center}
\end{table}

These are basic production cost parameters restricted by technological advance and complicated operations, and we do not expect them to change over short periods of time for all producers.
For the purpose of forecasting current year production, these parameters are revised monthly taken based on the market share of the month before.

The stochastic parameters are the risk-adjusted spot price $\rho_i(\xi)=\alpha(\xi)-\beta_i(\xi)$ and $H(\xi)={\rm diag}(h_i(\xi)+\gamma(\xi))$, where $\alpha(\xi)$ is benchmark price,  $\gamma(\xi)$ is stochastic supply discount, $\alpha_i$ and $\beta_i$ are the supply cost coefficients.

 For our experiments, we randomly choose $\zeta\in [0.05,0.1]$, and let $h_i=\beta_i=\zeta\times a_i$
 representing $5\%$ to 10\% of the unit production cost.
The data (ii) gives the crude oil price change due to different factors of contributions, namely contribution of demand $\Delta D$, supply $\Delta S$ and the residual contribution $\Delta R$.
Then, price change $\Delta {\rm price}$ is computed as follows:
\[
\Delta {\rm price}=\Delta D+\Delta S+\Delta R.
\]
%where $\Delta D,$ $\Delta S$ and $\Delta R$ can be positive or negative.
These contributions $\Delta D,$ $\Delta S$ and $\Delta R$ over certain period of time are uncertain. We assume that it can be described by random variable $\xi$ with unknown distribution, written as $d(\xi)$, $s(\xi)$ and $r(\xi)$.
For the purpose of our numerical tests, we formulate empirical distributions of historical data and use them as an approximation to the unknown distributions of different factors of contribution respectively.
Recall that in our model  the price is given by
\[
p(y(\xi), \xi)=\alpha(\xi) -\gamma(\xi) \sum^J_{j=1} y_j(\xi),
%p^k(\xi_\ell)-\gamma^k(\xi_\ell)T^{k}(y_{\xi_\ell}).
\]
in which the demand $D(\xi)$ is ignored as it would be a constant term and has no effect on solution.
Then, for any realization of $\alpha^k(\xi_\ell)$ of $k$-th day, it corresponds to
\[
\alpha^k(\xi_\ell)=\alpha^k_0(\xi_\ell)(1+d^k(\xi_\ell)+r^k(\xi_\ell)),
\]
where $\alpha^k_0$ is  the known price of prior day given in (iii),
 $d^k(\xi_\ell)$ and $r^k(\xi_\ell)$ are random scenarios taken from empirical distributions of $d(\xi)$ and $r(\xi)$, respectively.

%Supply $\sum_{i=1}^J y_i(\xi_\ell)$ which are taken from a uniformly distributed interval between $99{}$ and $101{}$ of yearly based daily production data from data (i).
It follows that, we can generate a set of data of
stochastic supply discount $\gamma(\xi_\ell)$
\[
%\gamma(\xi_\ell)=\frac{|\alpha^k(\xi_\ell)-\alpha^k_0|}{\sum_{i=1}^J y_i(\xi_\ell)},
\gamma(\xi_\ell)=\frac{|\alpha^k(\xi_\ell)-\alpha^k_0|}{\xi_\ell\bar{\eta}},
\]
where absolute value $|\cdot|$ ensures that increase in quantity has a negative influence on price,  $\xi_\ell \in [0.99,1.01]$ is uniformly distributed   and $\bar{\eta}$ is the total supply obtained from data (i).  We chose sample size $\nu=800$ of random variable $\xi$ for  both in-sample and  out-of-sample in the numerical simulation.

\begin{figure}
\centering
  % Requires
  \includegraphics[width=15cm]{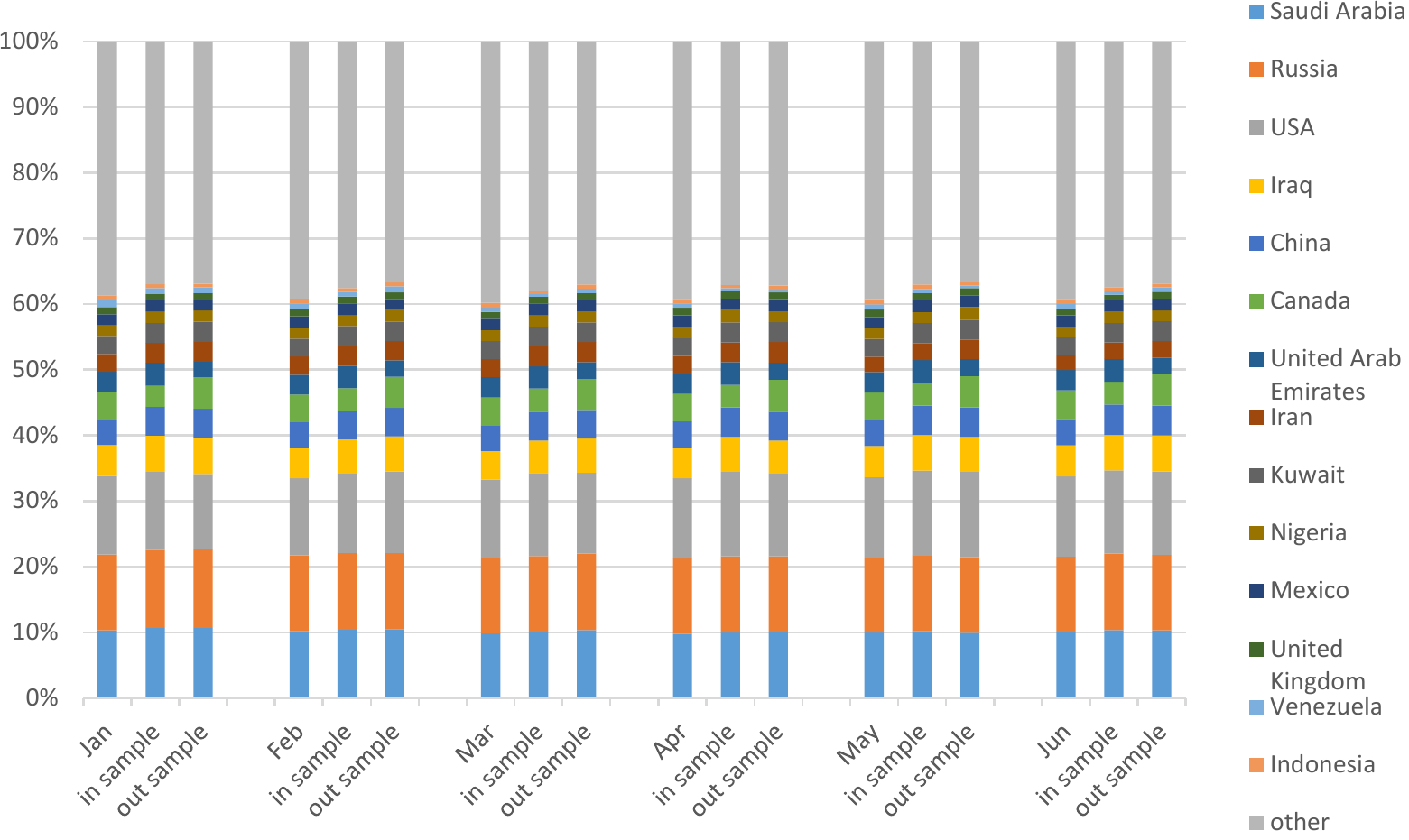}\\
    \caption{ Real, in sample and out sample monthly market shares of Jan to Jun in 2019
}\label{fig:short-term-2019-1}
    \includegraphics[width=15cm]{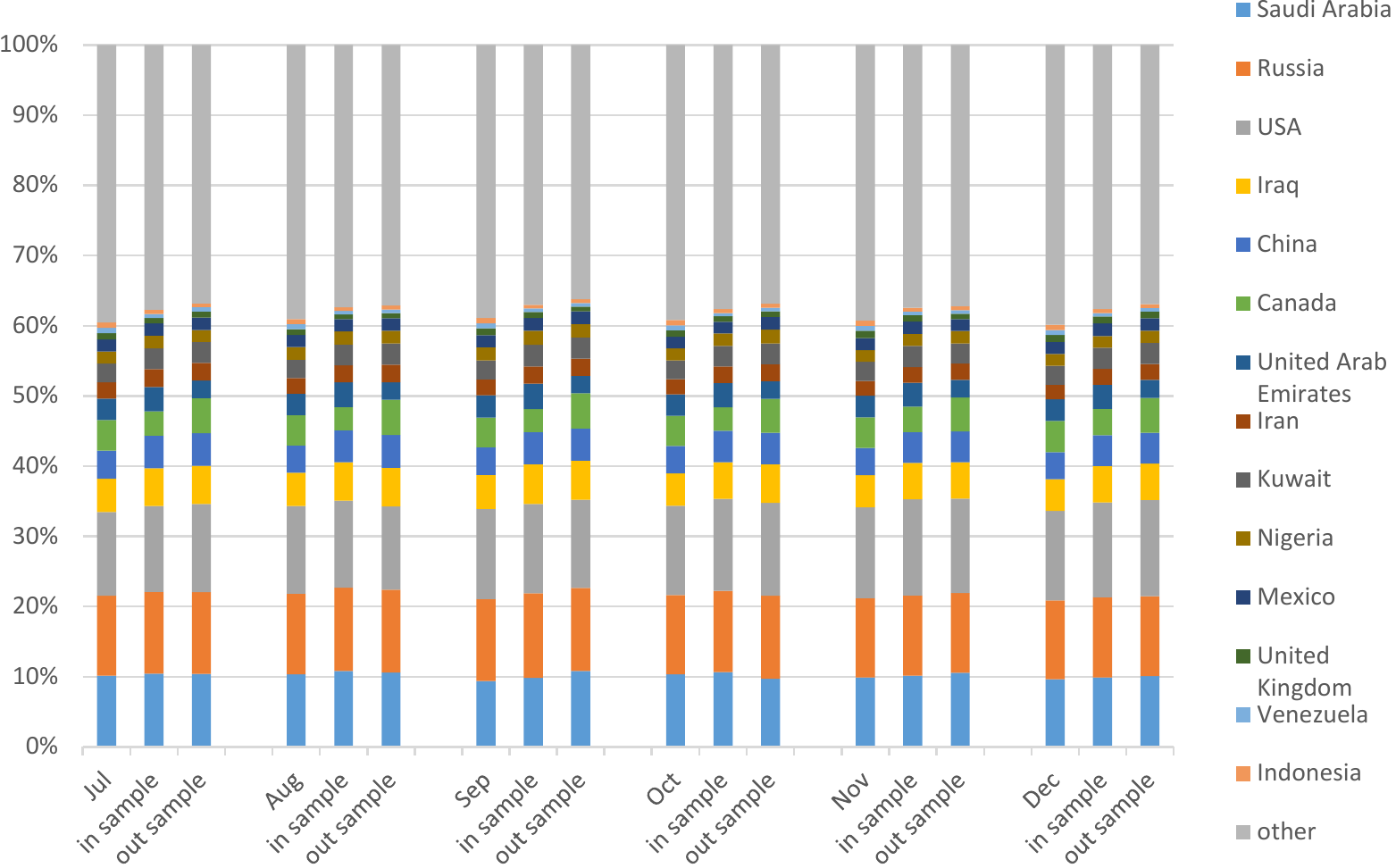}\\
  \caption{ Real, in sample and out sample monthly market shares of Jul to Dec in 2019
}\label{fig:short-term-2019-2}
\end{figure}

For long-term prediction (yearly market shares prediction), we refer interested readers to \cite{JSWC18} for more details. Here, we focus on short-term prediction, namely the monthly in-sample and out-of-sample market shares.
Table \ref{table:share-2019} gives average of daily market shares of producers in each month of 2019.
 Figures \ref{fig:short-term-2019-1} and \ref{fig:short-term-2019-2} display results for the recovered monthly market shares in 2019.
For each month, the first column is the real market share, while the second and third column are the in sample and out sample recovered results, respectively.
They show that our two-stage stochastic LCP model recovers and predicts the short-term real market shares from January 2019 to May 2020 very well.  Although global oil demand has been hit hard by COVID-19 and oil price has fell to historically low, our results show that a Nash-equilibrium for the global oil market share during the COVID-19 pandemic can be expected.

\begin{table}
\begin{center}
\begin{tiny}
\begin{tabular}{|c|c|c|c|c|c|c|c|c|c|c|c|c|}
  \hline
	& Jan &	Feb &	Mar	& Apr &	May	& Jun	& Jul	& Aug &	Sep &	Oct &	Nov &	Dec\\
\hline
Saudi Arabia&	10.31{}	&10.22{}	&	9.82{}	&	9.77{}	&	10.00{}	&	10.10{}	&	10.12{}	&	10.34{}	&	9.39{}	&	10.36{}	&	9.91{}	&	9.64{}	\\
Russia&	11.54{}	&	11.52{}	&	11.46{}	&	11.45{}	&	11.31{}	&	11.42{}	&	11.38{}	&	11.47{}	&	11.67{}	&	11.29{}	&	11.26{}	&	11.23{}	\\
USA	&11.95{}	&	11.75{}	&	11.96{}	&	12.27{}	&	12.33{}	&	12.25{}	&	11.98{}	&	12.49{}	&	12.83{}	&	12.74{}	&	12.93{}	&	12.76{}	\\
Iraq	&4.72{}	&	4.61{}	&	4.37{}	&	4.68{}	&	4.75{}	&	4.72{}	&	4.70{}	&	4.73{}	&	4.84{}	&	4.59{}	&	4.58{}	&	4.50{}	\\
China&	3.86{}	&	3.91{}	&	3.88{}	&	3.93{}	&	3.95{}	&	4.04{}	&	4.05{}	&	3.90{}	&	3.94{}	&	3.90{}	&	3.89{}	&	3.85{}	\\
Canada&	4.20{}	&	4.19{}	&	4.29{}	&	4.21{}	&	4.18{}	&	4.33{}	&	4.31{}	&	4.30{}	&	4.26{}	&	4.27{}	&	4.39{}	&	4.49{}	\\
UAE &	3.09{}	&	3.08{}	&	3.06{}	&	3.09{}	&	3.11{}	&	3.09{}	&	3.11{}	&	3.09{}	&	3.16{}	&	3.09{}	&	3.07{}	&	3.04{}	\\
Iran	&2.71{}	&	2.77{}	&	2.80{}	&	2.68{}	&	2.34{}	&	2.29{}	&	2.28{}	&	2.22{}	&	2.24{}	&	2.17{}	&	2.11{}	&	2.11{}	\\
Kuwait	&2.73{}	&	2.73{}	&	2.73{}	&	2.72{}	&	2.75{}	&	2.68{}	&	2.68{}	&	2.63{}	&	2.73{}	&	2.65{}	&	2.71{}	&	2.71{}	\\
Nigeria&	1.70{}	&	1.62{}	&	1.66{}	&	1.80{}	&	1.60{}	&	1.65{}	&	1.71{}	&	1.81{}	&	1.86{}	&	1.70{}	&	1.68{}	&	1.64{}	\\
Mexico	&1.63{}	&	1.72{}	&	1.70{}	&	1.70{}	&	1.69{}	&	1.70{}	&	1.69{}	&	1.70{}	&	1.75{}	&	1.66{}	&	1.71{}	&	1.70{}	\\
UK&	1.08{}	&	1.13{}	&	1.12{}	&	1.13{}	&	1.13{}	&	0.97{}	&	0.91{}	&	0.83{}	&	0.97{}	&	0.93{}	&	1.03{}	&	1.04{}	\\
Venezuela	&1.04{}	&	0.86{}	&	0.60{}	&	0.66{}	&	0.82{}	&	0.79{}	&	0.76{}	&  0.74{}	&0.73{}	&	0.70{}	&	0.70{}	&	0.71{}	\\
Indonesia	&0.77{}	&	0.76{}	&	0.76{}	&	0.71{}	&	0.77{}	&	0.73{}	&	0.75{}	&	0.74{}	&	0.75{}	&	0.74{}	&	0.73{}	&	0.73{}	\\
other	&38.65{}	&	39.15{}	&	39.79{}	&	39.19{}	&	39.27{}	&	39.24{}	&	39.56{}	&	39.03{}	&	38.89{}	&	39.21{}	&	39.29{}	&	39.87{}	\\
  \hline
\end{tabular}
\end{tiny}
 \caption{Average of daily market shares of producers in each month of 2019}\label{table:share-2019}
 \end{center}
\end{table}

\begin{table}
\begin{center}
\begin{tiny}
\begin{tabular}{|c|c|c|c|c|c|}
  \hline
	& Jan &	Feb &	Mar	& Apr &	May	\\
\hline
Saudi Arabia&	9.72{}	&	9.75{}	&	10.23{}	&	11.57{}	&	9.44{}	\\
Russia	&11.26{}	&	11.30{}	&	11.33{}	&	11.42{}	&	10.44{}	\\
USA	&12.72{}	&	12.85{}	&	12.76{}	&	12.28{}	&	11.23{}	\\
Iraq	&4.25{}	&	4.53{}	&	4.54{}	&	4.49{}	&	4.59{}	\\
China&	3.88{}	&	3.87{}	&	3.92{}	&	3.90{}	&	4.32{}	\\
Canada	&4.36{}	&	4.41{}	&	4.42{}	&	3.74{}	&	3.66{}	\\
UAE&	2.98{}	&	2.99{}	&	3.54{}	&	3.88{}	&	2.72{}	\\
Iran	&2.10{}	&	2.05{}	&	2.01{}	&	1.96{}	&	2.17{}	\\
Kuwait	&2.66{}	&	2.66{}	&	2.90{}	&	3.13{}	&	2.42{}	\\
Nigeria&	1.61{}	&	1.67{}	&	1.88{}	&	1.75{}	&	1.65{}	\\
Mexico	&1.72{}	&	1.75{}	&	1.77{}	&	1.75{}	&	1.83{}	\\
UK	& 0.97{}	&	0.95{}	&	1.02{}	&	1.00{}	&	1.05{}	\\
Venezuela	&0.73{}	&	0.75{}	&	0.73{}	&	0.71{}	&	0.76{}	\\
Indonesia	&0.73{}	&	0.72{}	&	0.72{}	&	0.72{}	&	0.79{}	\\
other	&40.32{}	&	39.75{}	&	38.24{}	&	37.68{}	&	42.93{}	\\
  \hline
\end{tabular}
\end{tiny}
 \caption{Average of daily market shares of producers in each month of January-May 2020}\label{table:share-2020}
 \end{center}
\end{table}
%\begin{figure}
%\centering
%  % Requires
%  \includegraphics[width=15cm]{figure/share_monthly_2019.pdf}\\
%    \caption{ Real monthly market shares in 2019
%}\label{fig:share-2019}
%\end{figure}

\begin{figure}
 \includegraphics[width=15cm]{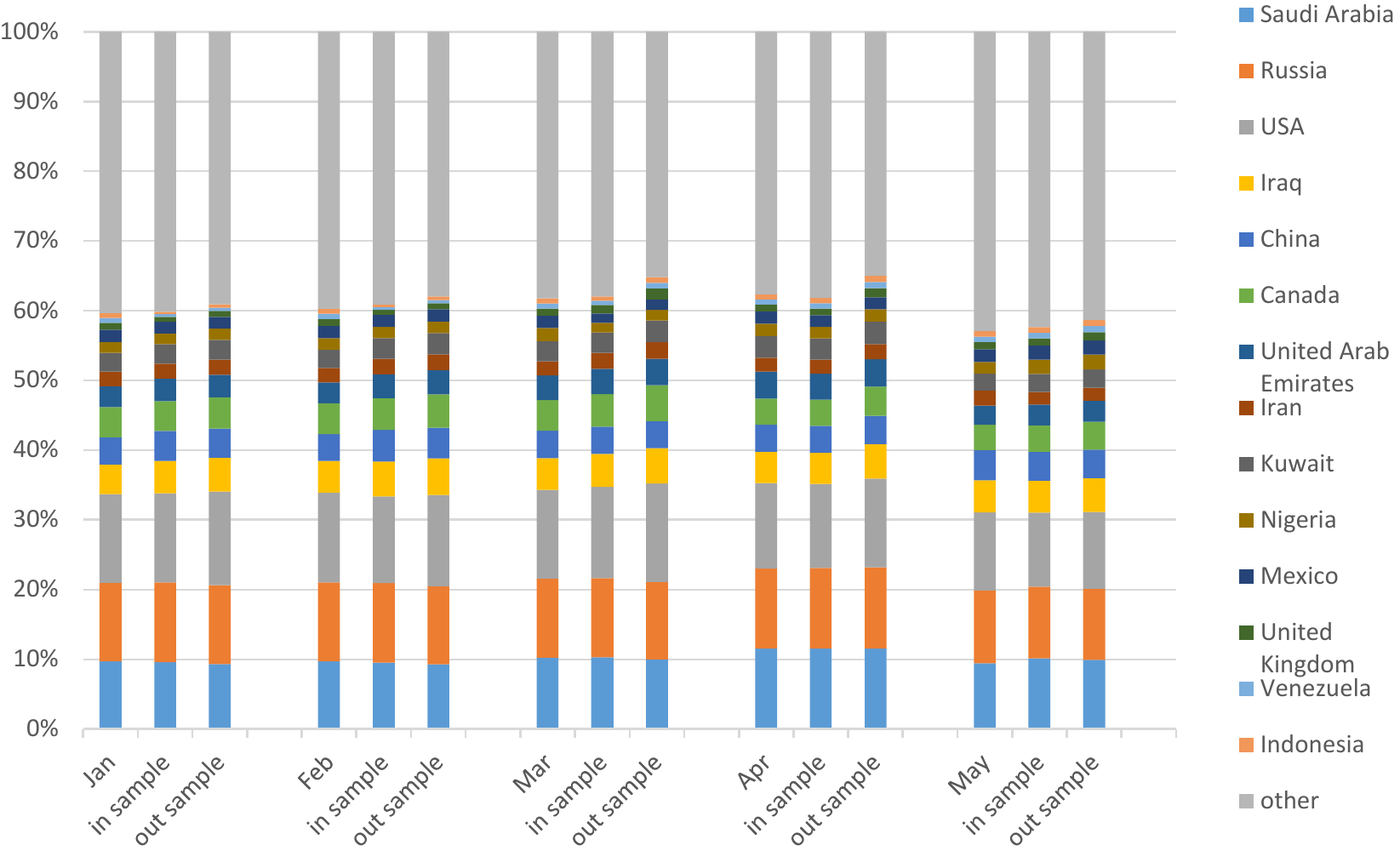}\\
  \caption{ Real, in sample and out sample monthly market shares of Jan to May in 2020
}\label{fig:short-term-2020}
\end{figure}

Of particular interests are what had happened in March and April among Russia, Saudi Arabia and U.S.A.
In particular, the Brent spot price fell from around \$70 per barrel to about \$50 since the identification of the pandemic and was believed to decrease further.
If no actions of change were taken by the major producers, what could have happened is that the high-cost producers would be forced to cut production because they are more vulnerable in low price environment.
On one hand, if the production cut is significant, the price will boost and the non-limited producers can be more profitable.
On the other hand, if the oil price continues to fall, the optimal production decision has little to change for the traditional producers with low unit cost, e.g., Russian, Saudi Arabia.
In indeed, with our model of two-stage stochastic game, the estimated market share given no strategic changes can be seen in accordance with the description above.
It is, given the estimated market share, the most rational decision for Russia is to refuse the production cut agreement.
Saudi Arabia, who believed to have the lowest unit production cost, followed the strategy of Russia immediately by offering price discount and increased its production.
For both countries, the low price environment has little effects in the sense of maintaining their market share respectively.
The same cannot be said for U.S.A., who has high cost in unit oil production, and is estimated to loss its market share if not to change its production strategy, presented by choosing non-zero $r$ in our model. It is most apparent from Table 2, that all the major producers responded in choosing their strategies to increase production quantities. The reality was counter intuitive at first glance, but those were in fact all rational decisions  and can be forecasted by our model

\section{Conclusion}
In this paper, we  model  the oil market share using the two-stage stochastic LCP (\ref{eq:LCP}) via the two-stage stochastic games (\ref{eq2})-(\ref{eq1}) with quadratic concave utility functions (\ref{qp1})-(\ref{qp3}).
  We show the existence and uniqueness of a Nash-equilibrium for the oil market share by the solution of the two-stage stochastic LCP (\ref{eq:LCP}).  Moreover, we propose the Alternating Block Algorithm (ABA) to find a solution of the
  two-stage stochastic LCP (\ref{eq:LCP}) with finite realizations. We derive the convergence theorems of the ABA and the PHA for
  solving (\ref{eq:LCP}) with finite realizations and show the ABA is much faster than the PHA by randomly generated problems.  We apply the new theoretical results and the ABA  to
the analysis of the oil market share from January 2019 to May 2020. Simulation results show that our model is effective for forecasting the oil market share under an uncertain environment of COVID-19  pandemic.  \\

{\bf Acknowledgements} 

Xiaojun Chen would like to thank Hong Kong Research Grant Council for grant PolyU153001/18P.\\

\vspace{0.1in}

{\bf Appendix,  Proofs of Lemma 3.1, Theorems 3.1, 4.1, 4.2}
\begin{appendix}
{\bf Proof of Lemma 3.1}
\begin{proof}
It is easy to see that $v(\xi)=(x, |\rho(\xi)|)^T\ge 0$ and $M(\xi)v(\xi)+q(x,\xi)\ge 0$, that is, $(x, |\rho(\xi)|)^T$
is a feasible solution of the LCP$(M(\xi),q(x,\xi))$. Hence from the feasibility and the positive semi-definiteness  of  $M(\xi)$, the LCP$(M(\xi),q(x,\xi))$ has at least one solution \cite[Theorem 3.1.2]{CoPaSt92}.

The LCP$(M(\xi),q(x,\xi))$ is the first order optimality condition of the strongly convex quadratic
program
%\begin{equation}\label{qp}
$$\begin{array}{cc}
\min & \frac{1}{2} y(\xi)^T(H(\xi) + \gamma(\xi)ee^T)y(\xi) +\rho(\xi)^Ty(\xi)\\
{\rm s.t.} &  0\le y(\xi) \le x.
\end{array}
\eqno{(\mbox{A.1})}
$$%\end{equation}
Because of the strong convexity, the first order optimality condition is necessary and sufficient for
the unique optimal solution  $y^*(\xi)$ of problem (A.1),  which is a fixed point of the fixed point problem
$$
y(\xi)=\Pi_{[0, x]}(y(\xi)-\rho(\xi)-  (H(\xi)+\gamma(\xi)ee^T)y(\xi)).
$$
Hence, $y^*(\xi)$ has the form of (\ref{fixed}).

By the definition of the projection, we can easily obtain
$$(\rho(\xi)+ (H(\xi)+\gamma(\xi)ee^T)y^*(\xi))_i$$
$$\left\{\begin{array}{ll}
=0 & \quad {\rm if} \, (y^*(\xi)-\rho(\xi)-(H(\xi)+\gamma(\xi)ee^T)y^*(\xi))_i \in (0, x_i)\\
\le 0 & \quad {\rm if} \, (y^*(\xi)-\rho(\xi)-(H(\xi)+ \gamma(\xi) ee^T)y^*(\xi))_i \ge x_i\\
\ge 0 & \quad {\rm otherwise.}
\end{array}
\right.
$$
Hence, the multiplier  $s^*(\xi)$ in (\ref{sols}) with $y^*(\xi)$ in (\ref{fixed})  is  a solution of the LCP$(M(\xi),q(x,\xi))$,
where for the case $y_i^*(\xi)<x_i$ it is from $s_i^*(\xi)= (\rho(\xi)+ (H(\xi)+\gamma(\xi)ee^T)y^*(\xi))_i=0$ and $s^*_i(\xi)(x_i-y^*_i(\xi))=0$, and for the case $y_i^*(\xi)=x_i$, it is from $y_i^*(\xi)(s^*_i(\xi)+(\rho(\xi)+ (H(\xi)+\gamma(\xi)ee^T)y^*(\xi))_i)=0.$

If $y_i^*(\xi)=x_i>0$, $s_i^*(\xi)$ in (\ref{sols}) is uniquely defined. If $y_i^*(\xi)=x_i=0$, then
$s_i^*(\xi)=\max(0, -(\rho(\xi)+ (H(\xi)+\gamma(\xi)ee^T)y^*(\xi))_i)=0$. Hence $v^*(\xi)=(y^*(\xi),s^*(\xi))^T
$ defined in (\ref{fixed})-(\ref{sols}) is the least norm solution of the  LCP$(M(\xi),q(x,\xi))$ and
the unique solution of LCP$(M(\xi),q(x,\xi))$ if $x>0$.
We complete the proof.
\end{proof}

{\bf Proof of Theorem 3.1}

\begin{proof}
We first prove (\ref{error_bound}).  Let
$$w(\xi)= (H(\xi)+\gamma(\xi)ee^T)y^*(\xi) + \rho(\xi) \quad {\rm and} \quad \bar{w}=
 (\bar{H}+\bar{\gamma}ee^T)\bar{u} + \bar{\rho}.$$
Then we have
$$y^*(\xi)=\Pi_{[0,x]}(y^*(\xi)-w(\xi)) \quad {\rm and} \quad \bar{u}=\Pi_{[0,x]}(\bar{u}-\bar{w}),$$
which implies
$${\rm mid} (y^*(\xi), y^*(\xi)-x, w(\xi))=0  \quad {\rm and} \quad {\rm mid} (\bar{u}, \bar{u}-x, \bar{w})=0,$$
where ``mid'' is the componentwise median operator. Following the proof of Lemma 2.1 in \cite{Chen_Wang}, there is a diagonal matrix $D=$diag$(d_1,\ldots,d_J)$ with $0\le d_i\le 1$ such that
$$0=(I-D)(y^*(\xi)-\bar{u})+D(w(\xi)-\bar{w}).$$
Hence we obtain
%\begin{equation}\label{Gamma}
$$(I-D+D(H(\xi)+\gamma(\xi)ee^T))(y^*(\xi)-\bar{u})=-D(\rho(\xi)-\bar{\rho} +(H(\xi)-\bar{H}+ (\gamma(\xi)-\bar{\gamma})ee^T)\bar{u}). \eqno{(\mbox{A.2})}$$
%\end{equation}
Since $H(\xi)+\gamma(\xi)ee^T$ is symmetric positive definite, by Theorem 2.7 in \cite{Chen_Xiang2}, we have
\begin{eqnarray*}
\max_{d\in[0,1]^J}\|(I-D+ D(H(\xi)+\gamma(\xi)ee^T))^{-1}D\|&=&\|(H(\xi)+\gamma(\xi)ee^T))^{-1}\|\\
&\le& \gamma(\xi)^{-1}\|(I+ee^T)^{-1}\|.
\end{eqnarray*}
Therefore, using $\|(I+ee^T)^{-1}\|=1$, $\|ee^T\|=J$, $\gamma(\xi)^{-1}\le \Gamma$ and $0\le \bar{u}\le x,$  we  obtain (\ref{error_bound}) from (A.2).

Next we show the continuity of the last $J$-components of the least norm solution of the LCP$(M(\xi), q(x,\xi))$.
Without loss of generality, assume that  $\sigma=\{i \in\{1,\ldots, J\} \, | \,  y_i^*(\xi)<x_i\} \neq \emptyset.$
From (\ref{error_bound}), for any $0<\epsilon <\min_{i\in \sigma}x_i-y^*_i(\xi)$, there is $\delta>0$ such that if $\|\rho(\xi)-\bar{\rho}\|+\|H(\xi)-\bar{H}\|+\|\gamma(\xi)-\bar{\gamma}\|<\delta$, then
$\|y^*(\xi)-\bar{u}\|<\epsilon.$ This implies for $i\in \sigma$,
$x_i-\bar{u}_i\ge x_i- y_i^*(\xi)-| y_i^*(\xi)-\bar{u}_i|>0.$ Hence, we have $\bar{t}_i=s_i^*(\xi)=0$.

For $i\not\in\sigma,$ that, is, $ y_i^*(\xi)=x_i$, we have $s^*_i(\xi)=-(\rho(\xi)+(H(\xi)+\gamma(\xi)ee^T)y^*(\xi))_i\ge 0$.
If $\bar{t}_i=0$, then $-(\bar{\rho}+(\bar{H}+\bar{\gamma} ee^T)\bar{u})_i\le 0$, otherwise
$\bar{t}_i=-(\bar{\rho}+ (\bar{H}+  \bar{\gamma} ee^T)\bar{u})_i\ge 0$. Hence, using  (\ref{error_bound}), we have
\begin{eqnarray*}
&&|s^*_i(\xi)-\bar{t}_i|\\
&\le&|(\rho(\xi)+(H(\xi)+\gamma(\xi)ee^T)y^*(\xi))_i -(\bar{\rho}+  (\bar{H}+\bar{\gamma} ee^T)\bar{u})_i|\\
&\le&\|\rho(\xi)-\bar{\rho}\|+(\|H(\xi)\|+\gamma(\xi)J)\|y^*(\xi)-\bar{u}\|
+\|x\|(J\|\gamma(\xi)-\bar{\gamma}\|+\|H(\xi)-\bar{H}\|)\\
&\le& (L+(\|H(\xi)\|+\gamma(\xi)J)\Gamma)(\|\rho(\xi)-\bar{\rho}\|+\|H(\xi)-\bar{H}\|+\|\gamma(\xi)-\bar{\gamma}\|)
),
\end{eqnarray*}
where $L\ge \max\{1, J\|x\|\}.$
Hence, the solution of the LCP$(M(\xi), q(x,\xi))$ is continuous with respect to $H(\xi)$, $\rho(\xi)$ and $\gamma(\xi)$. We complete the proof.

\end{proof}

{\bf Proof of Theorem 4.1}

\begin{proof}
Let $\Lambda=$diag($\nu, 1,\ldots, 1)\in \mathbb{R}^{n\times n }$.  It is ease to verify that $\mathbf{v}^*$ is a solution of
the  LCP$(\Lambda\mathbf{M}$,$\Lambda \mathbf{q})$ if and only if $\mathbf{v}^*$ is a solution of  the  LCP$(\mathbf{M}$,$\mathbf{q})$.

Since
$z^T(\Lambda \mathbf{M}+ (\Lambda \mathbf{M})^T)z\ge 0$ for any $z\in \mathbb{R}^{n}$, the matrix $\Lambda\mathbf{M}$ is positive semi-definite.  Moreover,
$\mathbf{v}=(x, x, |\rho(\xi_1)|,\ldots, x, |\rho(\xi_\nu)|)^T$
with $ x_i > \max\{0,  (\frac{1}{\nu}\sum^\nu_{\ell=1} |\rho(\xi_\ell)| -a )_i/c_i\}$, $i=1,\ldots, J$ is a feasible solution of
the  LCP$(\Lambda\mathbf{M}$,$\Lambda \mathbf{q})$. Hence the feasibility implies that  LCP$(\Lambda\mathbf{M}$,$\Lambda \mathbf{q})$ has at least one solution $\mathbf{v}^*$ \cite{CoPaSt92}.
If $\mathbf{v}^*$ has the component $x^*>0$, then by Lemma 3.1, the component $(v^*(\xi_\ell),\ldots,
v^*(\xi_\nu))$ of $\mathbf{v}^*$ is uniquely dependent on $x^*>0$. Moreover,
we have
$$(C+re^T)x^* + \frac{1}{\nu}\sum^\nu_{i=\ell}D_\ell(\rho(\xi_\ell)+ (H(\xi_\ell)+\gamma(\xi_\ell) ee^T)x^*)+a=0,$$
where $D_\ell$ is a diagonal matrix with diagonal elements being 0 or 1. Since $C+re^T$ and $H(\xi_\ell)+\gamma(\xi_\ell) ee^T$ are positive definite, $x^*$ is the unique solution of this system of equations.  Hence,
the LCP$(\Lambda \mathbf{M}$,$\Lambda \mathbf{q})$ has  at most one solution with $x^*>0$.

For any $z\in \mathbb{R}^{n}$, we can easily find
$$\frac{1}{\nu}
z^T(\Lambda \mathbf{M}+ (\Lambda \mathbf{M})^T)z=\frac{1}{\nu} \sum_{\ell=1}^\nu z_\ell^T\left(\begin{array}{cc}
C+re^T&   -B \\
B^T & M(\xi_\ell)\end{array}\right)z_\ell^T  \ge 0,$$
where $z_\ell=(x, v(\xi_\ell))^T$ for $\ell=1,\ldots, \nu$. Moreover, the matrix $ \left(\begin{array}{cc}
C+re^T &   -B \\
B^T & M(\xi_\ell)\end{array}\right)$ is positive semi-definite for any $\xi_\ell$. Hence from \cite{RoSun18a}, the sequence $\{z^k\}$ generated by Algorithm 1 converges to a solution of the LCP$(\Lambda \mathbf{M}$,$\Lambda \mathbf{q})$.   We complete the proof.
\end{proof}

{\bf Proof of Theorem 4.2}

\begin{proof} From Step 2 of Algorithm 2 and the definition of a solution, we have
$$0 \le  x^{k+1} \bot  (C+re^T)x^{k+1} -\frac{1}{\nu} \sum^\nu_{\ell=1} s^k_\ell + a \ge 0 \quad {\rm and} \quad  0 \le  x^* \bot  (C+re^T)x^* -\frac{1}{\nu} \sum^\nu_{\ell=1} s^*_\ell + a \ge 0,$$
which imply
\begin{eqnarray*}
0 &\ge &(x^{k+1}-x^*)^T((C+re^T)(x^{k+1}-x^*)+ \frac{1}{\nu} \sum^\nu_{\ell=1} (s^k_\ell-s_\ell^*)\\
&\ge & \|C_2(x^{k+1}-x^*)\|^2 - \|C_2(x^{k+1}-x^*)\|\|C_2^{-1} \|\frac{1}{\nu} \sum^\nu_{\ell=1} \|s^k_\ell-s_\ell^*\|.
  \end{eqnarray*}
Let $\mu>0$ such that $\Omega:=\{ x \, | \, \|C_2(x-x^*)\|\le \mu\} \subset \mathbb{R}_{++}^J$.  By Theorem 2.1 in \cite{Chen_Xiang3} and Lemma 3.1, for any $x^k \in \Omega$, the LCP$(M(\xi_\ell), q(x^k,\xi_\ell))$ has a unique solution $v^k_\ell=(y^k_\ell, s^k_\ell)^T$, which is a Lipschitz continuous function  with  the Lipschitz constant $\sigma_\ell$ of $x$. Hence, from $q(x,\xi)=(\rho(\xi),x)^T$, we have
$$\|s^k_\ell-s_\ell^*\|\le \|v^k_\ell-v^*\|\le \sigma_\ell \|q(x^k,\xi_\ell)-q(x^*,\xi_\ell)\|=\sigma_\ell \|x^k-x^*\|$$
and
\begin{eqnarray*}
\|C_2(x^{k+1}-x^*)\|&\le& \|C_2^{-1}\| \frac{1}{\nu} \sum^\nu_{\ell=1} \|s^k_\ell-s_\ell^*\|
\le  \|C_2^{-1}\| \frac{1}{\nu} \sum^\nu_{\ell=1} \sigma_\ell \|x^k-x^*\|\\
&\le& \|C_2^{-1}\|^2 \sigma \|C_2(x^k-x^*)\|.
\end{eqnarray*}
Since $\|C_2^{-1}\|^2 \sigma<1$,  we  obtain $x^{k+1} \in \Omega$, and the convergence of $\{x^k\}$ to
$x^*$. Using Lemma 3.1 again, the sequence $(x^k, y_1^k, s^k_1, \ldots, y^k_\nu, s^k_\nu)^T $ generated by the ABA converges to the solution of the LCP$(\mathbf{M}$,$\mathbf{q})$.
\end{proof}

\end{appendix}

\end{document}